\documentclass[oneside,12pt,reqno]{amsart}
\usepackage{amssymb, amscd}
\usepackage[all]{xy}
\def\BZ{{\mathbb{Z}}}

\def\ll{{\mbox{\boldmath $\ell$}}}

\newcommand{\hh}{{\mbox{\boldmath $h$}}}

\def\NN{{\mathcal N}}

\def\EE{{\mathcal E}}

\begin{document}
\begin{titlepage}
\vspace{1.cm}
\begin{flushright}
LMU-ASC 57/05\\
\end{flushright}
\vspace{1cm}

\begin{centering}
\vspace{1cm} {\Large {\bf Crossed Module Bundle Gerbes; Classification, String
Group and Differential Geometry}}\\
\vspace{1cm}

{\bf Branislav Jur\v co} \\ \vspace{0.1cm} Theoretical Physics, LMU Munich\\
Theresienstr. 37, 80333 Munich\\ Germany \vspace{.25cm}

\vspace{1.5cm}

{\bf Abstract}
\\
\end{centering}
\vspace{0.5cm} 
\noindent We discuss nonabelian bundle gerbes and
their differential geometry using simplicial methods. Associated to
any crossed module $(H\to D)$ there is a simplicial group $N\mathcal
C_{(H\to D)}$, the nerve of the 1-category defined by the crossed
module and its geometric realization $|N\mathcal
C_{(H\to D)}|$. Equivalence classes of principal bundles with
structure group $|N\mathcal
C_{(H\to D)}|$ are shown to be one-to-one with stable
equivalence classes of what we call crossed module 
bundle gerbes. We can also associate to a crossed module a 2-category
$\tilde{\mathcal C}_{(H\to D)}$. Then there are two equivalent ways how to view 
classifying spaces of
$N\mathcal C_{(H\to D)}$ -bundles and hence of $|N\mathcal
C_{(H\to D)}|$-bundles and crossed module bundle gerbes. We can either apply the $W$-construction to 
$N\mathcal C_{(H\to D)}$ or take the nerve of the 2-category 
$\tilde{\mathcal C}_{(H\to D)}$. 
We discuss the string group and string structures from
this point of view. Also a simplicial principal bundle can be
equipped with a simplicial connection and a $B$-field. It is shown
how in the case of a simplicial principal $N\mathcal C_{(H\to
D)}$-bundle these simplicial objects give the bundle gerbe
connection and the bundle gerbe $B$-field.

\vspace{3cm}
\begin{flushleft}
{\ttfamily email: jurco@lmu.de}
\end{flushleft}
\end{titlepage}

\setcounter{page}{1}

\section*{0. Introduction}
Nonabelian gerbes arose in the realms of nonabelian cohomology \cite{Gi},
\cite{BGr} and
higher category \cite{Br}. Their differential geometry was described thoroughly
by Breen and Messing \cite{BM} from the algebraic geometry point of view
(see \cite{Attal} for the combinatorial description).
In \cite{ACJ} nonabelian bundle gerbes, generalizing the nice concept of an
abelian bundle gerbe \cite{Murray}, were introduced. These have to be shown
(along with their connections and curvings) very natural objects in classical
fibre bundle theory. There a is hope that in this form gerbes can be useful in
physics (see e.g. examples of higher Yang-Mills theories \cite{B} and anomaly
cancellation of M5-branes \cite{AJ}). Closely related to crossed modules bundle
gerbes are 2-bundles introduced in \cite{Bartels} and discussed together with their
connections and curvings in \cite{BS}.

In this paper we discuss classification of bundle gerbes associated with
crossed modules. These are bundle gerbes equipped with modules in the
terminology
of \cite{ACJ}. This is done using some well-known simplicial constructions.
Then the relation between simplicial principal bundles and crossed module bundle
gerbes is used to describe bundle gerbe connections and curvings 
in simplicial language.

The first section is devoted to simplicial principal bundles. We describe them
as twisted Cartesian products following \cite{May} and recall the construction
of the universal bundle.

Connections on simplicial bundles are introduced in Section 2. This is done in a
straightforward way, which we believe, is the relevant one for our purposes. Next
we shortly discuss the corresponding notion of a curvature.

Our task in Section 3 is to define the next in an infinite sequence 
of relevant differential geometric objects associated with simplicial principal
bundles, the $\bar B$-field.

In Section 4 we describe some simplicial constructions related to a crossed module
$(H\to D)$.
We can view a crossed module as a 1-category (actually 1-groupoid)
${\mathcal C}_{(H\to D)}$
or as a 2-category (actually a 2-groupoid) $\tilde {\mathcal C}_{(H\to D)}$.
We can form the corresponding nerves $N{\mathcal C}_{(H\to D)}$ and
$N\tilde {\mathcal C}_{(H\to D)}$ respectively.
The geometric realization $|N{\mathcal C}_{(H\to D)}|$ is the classifying space of
$H$-principal bundles with a chosen trivialization when we change the structure
group from $H$ to $D$. If $H$ and $D$ are Lie groups ${\mathcal C}_{(H\to D)}$
is a simplicial Lie group and its geometric realization
$|N{\mathcal C}_{(H\to D)}|$is a topological group. String
group of \cite{BCSS}, \cite{H} is an example. We remark on how the construction
of \cite{BCSS} relates to the one of Stolz and Teichner \cite{ST}.

Crossed module bundle gerbes are introduced in
Section 5. The geometric realization $|N\tilde {\mathcal C}_{(H\to D)}|$ of the nerve of the 2-category related to
our crossed module gives the
classifying space of such bundle gerbes. These are shown to be 
the same as principal
bundles with structure group $|N{\mathcal C}_{(H\to D)}|$. 
In particular string structures \cite{ST} can be described equivalently
in terms of nonabelian bundle gerbes.
Locally crossed module bundle gerbes can be described using simplicial maps 
between the nerve of the 
1-category coming from an open covering of the manifold and a nerve of the 
2-category associated with a crossed module (which is equal to the classifying
space $\overline{W} N{\mathcal C}_{(H\to D)}$ of 
principal $N{\mathcal C}_{(H\to D)}$-bundles). Here we have to mention closely
related work of D. Stevenson \cite{S}.
In the last section we describe how the connection and $\bar B$-field on
simplicial principal $N{\mathcal C}_{(H\to D)}$-bundle give rise
to a connection and an $B$-field on the corresponding crossed module bundle 
gerbe.
Here as well in sections 2 and 3 we work in the category of manifolds. However as J. Baez
pointed out it might be more
appropriate to work in the category of ``smooth spaces" studied in the appendix
of \cite{BS}.

Some generalization to the case of bigroupoids will be given in thesis of I.
Bakovi\'c \cite{igor}.

Finally we should mention \cite{BGr} again. Here the group 
$|N{\mathcal C}_{(H\to D)}|$ and classifying space 
$\overline{W} N{\mathcal C}_{(H\to D)}$ are discussed. I thank to D. Stevenson
for pointing out this to me. Also I am very much indebted A. Henriques for 
help with sections 4 and 5.

\section{Simplicial principal bundles}
We start by recalling some relevant properties of simplicial
principal bundles following mainly \cite{May}. Let $\pi:P\to X$ to
be a simplicial (left) principal $G$-bundle, with $P$ and $X$
simplicial sets  and $G$ a simplicial group. As usually we will use
$\partial_i$ and $s_i$ for the corresponding face and degeneracy
maps. In the rest of the paper we always assume, without spelling it
out explicitly, $P \to X$ to posses a pseudo-cross section $\sigma:
X \to P$ such that $\pi\sigma = id_X$ and $\partial_i \sigma =\sigma
\partial_i$ if $i>0$ and $s_i \sigma =\sigma s_i$ if $i\geq 0$.
Associated with a pseudo-cross section $\sigma$ we have the twisting
function $\tau: X_n \to G_{n-1}$ defined as
$$\partial_0 \sigma (x)=\tau(x).\sigma (\partial_0 x).$$

We will use the following description of $G$-bundles which we
alternatively can use as a definition.
\subsection{Twistings}
To make this section self-contained we have to describe the twisting first.
For a function $\tau: X_n \to G_{n-1}$ to be a twisting the following
conditions should be fulfilled:
$$ \partial_0\tau(x)= \tau(\partial_1 x)(\tau(\partial_0
x))^{-1},$$
$$\partial_i\tau(x)=\tau(\partial_{i+1} x)\hskip0.4cm {\rm for}
\hskip0.4cm i>0,$$
$$s_i\tau(x)= \tau(s_{i+1} x)\hskip0.4cm {\rm for}
\hskip0.4cm i\geq0,$$
$$ \tau(s_0 x)=e_n\hskip0.4cm {\rm for}
\hskip0.4cm x \in X_n.$$
\subsection{Principal bundles as twisted Cartesian products}
A principal $G$-bundle $p:P\to X$
with a pseudo-cross section can be identified with the simplicial
set $P(\tau)= G\times_{\tau}X$, which satisfies
$$ P(\tau)_n= G_n\times X_n$$ and has the following face and degeneracy maps
\\

$(i)$ $\hskip4.2cm \partial_i(g,x)= (\partial_i g,\partial_i x) \hskip0.4cm {\rm for}
\hskip0.4cm i>0,$

$(ii)$ $\hskip4.1cm \partial_0(g,x)=
(\partial_0 g .\tau(x),\partial_0 x), $

$(iii)$ $\hskip4cm s_i(g,x)= (s_i g, s_i x) \hskip0.4cm {\rm for}
\hskip0.4cm i\geq 0$.
\\

\noindent
Moreover there is a canonical choice for the pseudo-cross section
$\sigma(x)=(e_n, x)$, $x\in  X_n$ and $e_n$ the identity in $G_n$.

Equivalence of two $G$-bundles $P(\tau)$ and $P(\tau')$ over the
same $X$  is described in terms of twisting as follows.
\subsection{Equivalence of principal bundles}
We call two twistings $\tau'$ and $\tau$ equivalent if there exists
a simplicial map $\psi: X \to G$ such that
$$
\partial_0\psi(x).\tau'(x)= \tau(x).\psi(\partial_0x),$$
$$ \partial_i\psi(x)=\psi(\partial_i x) \hskip0.4cm {\rm if}\hskip0.4cm i>0,$$
$$s_i\psi(x)=\psi(s_i x) \hskip0.4cm {\rm if}\hskip0.4cm i\geq 0.$$

It will be convenient to introduce the equivariant map $\bar{\sigma}: P\to G$,
$ \bar{\sigma}(gp)=g.\bar{\sigma}(p)$ by the
equation $p=\bar\sigma(p)\sigma(x)$. In the rest we will always assume the
canonical choice of the pseudo-cross section is made in which case
$\bar\sigma(g_n,g_{n-1},\ldots, g_0)= (g_n)$.
We have
$$\partial_0\bar \sigma (p) = \bar \sigma(\partial_0 p)\tau(x)^{-1}.$$

As with ordinary bundles simplicial bundles can be pulled back and
their structure groups can be changed using simplicial group
homomorphisms. Pseudo-cross sections and twistings transform under
these operations in the usual way.

\subsection{Universal $G$-bundle}
There is a canonical choice of the classifying space of $G$-bundles denoted as
$\overline{W}G$ and constructed as follows. $\overline{W}G_0$ has one
element $\ast$ and $\overline{W}G_n=G_{n-1}\times G_{n-2}\times \ldots \times
G_0$ for $n>0$. Face and degeneracy maps are
$$s_0(\ast) = (e_0), \hskip0.4cm  \hskip0.4cm \partial_i(g_0)=\ast
\hskip0.4cm {\rm for}\,\,\,\,i=0\,\,{\rm or}\,\,1$$
and
$$\partial_0 (g_{n}, \ldots g_0) =(g_{n-1},\ldots g_0),$$
$$\partial_{i+1} (g_{n}, \ldots g_0) =(\partial_ig_{n},\ldots ,
\partial_1g_{n-i+1},\partial_0g_{n-i}.g_{n-i-1}, g_{n-i-2} ,\ldots,  g_0),$$
$$s_0(g_{n-1}, \ldots ,g_0)=(e_n, g_{n-1},\ldots, g_0),$$
$$s_{i+1}(g_{n-1}, \ldots ,g_0)=(s_ig_n,\ldots, s_0g_{n-i},e_{n-i},g_{n-i-1},
\ldots, g_0),$$
if $n>0$. With the choice of a twisting given by
$$\tau(g_{n-1}, \ldots ,g_0)=g_{n-1}$$ we have the universal $G$-principal
bundle
$$WG = G\times_\tau \overline{W}G.$$

As with ordinary bundles we have that $WG$ is contractible
and is universal in the following sense.

\subsection{Theorem} { \em Let us assign to any simplicial
map $$f: X \to \overline WG$$ the induced bundle $f^*(WG) \to X$.
This defines a one-to-one correspondence between homotopy classes of
maps $[X, \overline WG]$ and the equivalence classes of principal
$G$-bundles over the base $X$.}

\section{Simplicial connection, curvature}
Here we introduce the notion of a connection on a simplicial bundle. Of course
now we assume that $G$ is a simplicial Lie group and $P$ and $X$ are
simplicial manifolds. Also all maps an actions are smooth. We use the shorthand
notation $\Omega^k(Y) \otimes {\rm Lie}(G)$ for the collection of all ${\rm
Lie}(G_n)$-valued $k$-forms on $Y_n$ for all $n$ and any simplicial manifold $Y$.
Here of course ${\rm Lie}(G)$ is the corresponding simplicial Lie
algebra ${\rm Lie}(G)_n= {\rm Lie}(G_n)$ with the induced face and degeneracy
maps. For purposes of this paper the following definition of a simplicial
connection seems to be adequate.

\subsection{Definition} Let $\mathcal A \in \Omega^1(P) \otimes {\rm Lie}(G)$ be a
collection of one forms $\mathcal A_n \in \Omega^1(P_n) \otimes {\rm
Lie}(G_n)$. We call $\mathcal A$ a connection on the simplicial
principal $G$-bundle $P\to X$ if it fulfils the following
conditions:
\\

$(i)\hskip4cm \partial^*_i \mathcal A = \partial_i \mathcal A
\hskip0.4cm {\rm and}\hskip0.4cm s^*_i \mathcal A = s_i \mathcal A$
\\

\noindent where  $\partial^*_i$ on the left is the pullback of the
face map acting on the one-form part of $\mathcal A$ and $\partial_i
\mathcal A$ on the right is the simplicial Lie algebra face map
acting on the simplicial Lie algebra part of $\mathcal A$ and
similarly for degeneracies
\\

$(ii)$ $\mathcal A$ is equivariant with respect to the left
$G$-action on $P$
$$ g^* \mathcal A = g \mathcal A g^{-1}$$
and
\\

$(iii)$ its pullback to the fibre under $\bar \sigma: P \to G$
is the Cartan-Maurer form $g dg^{-1}$, i.e. the
collection of elements $ g_n dg_n^{-1} \in \Omega^1(G_n) \otimes {\rm
Lie}(G_n)$.
\subsection{Local connection forms}
Let us consider a collection of one forms
$A \in \Omega^1(X) \otimes {\rm Lie}(G)$ with the property
$$\partial_0 A =\tau\partial_0^*A\tau^{-1} + \tau d\tau^{-1},$$
$$\partial_i^*A=\partial_i A\hskip0.4cm {\rm for}\hskip0.4cm i>0$$
and
$$s_i^*A=s_i A\hskip0.4cm {\rm for}\hskip0.4cm i\geq0\,.$$
We call such an $A$ a local connection.

The following proposition is obvious.
\subsection{Proposition}
{\em Any connection $\mathcal A$ is of the form
$$\mathcal A = \bar \sigma A \bar \sigma^{-1} + \bar \sigma d \bar \sigma^{-1}$$
with
$$A = \sigma^* \mathcal A.$$}

Pullbacks and change of the structure group work as usually.
\subsection{Curvature}
Curvature is defined exactly the same way as in the case of ordinary
bundles. It is a collection of two forms $\mathcal F \in \Omega^2(P)
\otimes {\rm Lie}(G)$ defined as $\mathcal F = d\mathcal A +
\mathcal A\wedge \mathcal A$ and it has the following properties:
\\

$(i)\hskip4cm \partial^*_i \mathcal F = \partial_i \mathcal F
\hskip0.4cm {\rm and}\hskip0.4cm s^*_i \mathcal F= s_i \mathcal F$
\\

$(ii)$ $\mathcal F$ is equivariant with respect to the left
$G$-action on $P$
$$ g^* \mathcal F = g \mathcal F g^{-1}$$
and \\

$(iii)$ $\mathcal F$  is of the form $\mathcal F=\bar \sigma F
\bar \sigma^{-1}$ with $F \in \Omega^2(X) \otimes {\rm Lie}(G)$,
i.e. it is horizontal.

Of course $F=dA+ A\wedge A$.

Let us note that
$$\partial_0 F = \tau \partial_0^* F \tau^{-1} \hskip0.4cm {\rm
and}\hskip0.4cm \partial_iF =  \partial_i^* F \hskip0.4cm {\rm
for}\hskip0.4cm i>0$$ and
$$s_iF = s^*_i F \hskip0.4cm {\rm for}\hskip0.4cm i\geq0\,.$$
\section{$\bar B$-field}
Let
$$\bar G_0 = 0\,,$$
$$\bar G_n = \ker \partial_1\ldots \partial_{n}\subset G_n.$$
Let us note that $\partial_0\bar G_{n+1}$ is a normal subgroup in
$G_n$. Also $\partial_i \bar
G_{n+1} \subset \bar G_{n}$ for $i>0$. From now on we will assume that there
exists
an action of $G_n$ on $\bar
G_{n+1}$; $g_n \times \bar g_{n+1}
\mapsto \,^{g_n}\bar g_{n+1}$ such that
$$\partial_0(\,^{g_n}\bar g_{n+1})={g_n}\partial_0(\bar g_{n+1})g_n^{-1}$$
and
$$\,^{\partial_0 \bar g_{n+1}}\bar g_{n+1}'=\bar g_{n+1}\bar g_{n+1}'
\bar g_{n+1}^{-1}.$$
These conditions will be automatically satisfied in the next sections when we
consider simplicial groups with  simplicial homotopy groups $\pi_i (G) = 0$, for
$i\geq 2$.

\subsection{Definition} $\bar B$-field is a collection of
two-forms $\bar B_{n+1} \in\Omega^2(X_n) \otimes (\bar G_{n+1})$
such that
$$\,^\tau\partial^*_0\bar B =
\partial_1\bar B$$ and $$\partial^*_{i}\bar B =
\partial_{i+1} \bar B \hskip0.4cm {\rm for} \hskip0.4cm i>0$$
and
$$s^*_{i}\bar B =
s_{i+1} \bar B \hskip0.4cm {\rm for} \hskip0.4cm i\geq 0.$$

Finally we introduce collection of two forms $\nu \in \Omega^2(X)
\otimes {\rm Lie}(G)$ as
$$\nu_n = F_n + \partial_0 \bar B_{n+1}.$$
Obviously $\nu$ has the same properties with respect to face and
degeneration maps as $F$.
\subsection{Remark}
Of course there is no reason to stop with connection $A$ and $B$-field here.
One can introduce $C$-field etc. ad infinitum. We will however not to do so
here as we are really interested only in simplicial groups which are algebraic
models of homotopy 2-type (crossed modules). Also it seems that it would be more
proper to treat all this fields together, see remark 6.5.

\section{Crossed modules}

\subsection{Definition}
Let $H$ and $D$ be two Lie groups. We say that $H$ is a crossed $D$-module if there
is a group homomorphism $\alpha: H\to D$ and an action of $D$ on $H$ denoted as
$(d,h)\mapsto \,^d\hskip-0cm h$ such that
$$^{\alpha(h)}h' = hh'h^{-1}\hskip0.4cm {\rm{for}} \hskip0.4cm h,h'\in H$$
and
$$\alpha(^dh)= d \alpha(h)d^{-1} \hskip0.4cm {\rm{for}} \hskip0.4cm h\in H,
d\in D.$$
holds true.

We will use the  the following notation  $(H\to D)$ for a crossed
module. If the groups are infinite dimensional we have assume that
these are Frech\'et Lie Groups.

There are two canonical categorial construction associated with any crossed
module.

\subsection{Crossed module as a 1-category}
Let us denote ${\mathcal C}_{(H \to D)}$ the (topological) category
with objects being group elements $d\in D$ and morphisms (1-arrows)
group elements $(h,d)$ of the semidirect product $H\rtimes D$.

As with any category we can now form the simplicial space, the nerve
$N{\mathcal C}_{(H \to D)}$ of ${\mathcal C}_{(H \to D)}$ and its
geometric realization $|N{\mathcal C}_{(H\to D)}|$. The nerve is
naturally a simplicial Lie group and its geometric realization becomes
naturally a topological group \cite{BCSS}. We will use the following 
pictorial representation for the simplicial
group $N{\mathcal C}_{(H \to D)}$:
$$\stackrel{d_0}{\longrightarrow}$$ for the zeroth component, 
\[\xymatrix{&\ar@2[dd]^-{h_{01}}\\
\ar@/^1.1cm/[rr]^-{d_1}\ar@/_1.1cm/[rr]_-{d_0}&&\\
&&}\]
for the first component, 
\[\xymatrix{&&\ar@2[d]^-{h_{12}}\\
&\ar@/^1.1cm/[rr]^-{d_2}\ar[rr]^{\hskip-0.6cm{d_1}}
\ar@/_1.1cm/[rr]_-{d_0}&\ar@2[d]^-{h_{01}}&&
\\&&}
\] for the second component etc., with the obvious face and degeneracy maps.
%\[\xymatrix{&&&&&\ar@2[dd]^-{h_{01}}&&&&&\ar@2[d]^-{h_{12}}\\
%\ar[r]^-{d_0}& &
%&\ar@<2pt>[l]\ar@<-2pt>[l]&\ar@/^1.1cm/[rr]^-{d_1}
%\ar@/_1.1cm/[rr]_-{d_0}
%& & &
%&\ar@<4pt>[l]\ar@<-4pt>[l]\ar[l]&\ar@/^1.1cm/[rr]^-{d_2}\ar[rr]^{\hskip-0.6cm{d_1}}
%\ar@/_1.1cm/[rr]_-{d_0}&\ar@2[d]^-{h_{01}}&&\ldots
%\\&&&&&&&&&&&}
%\] 

The (opposite) group multiplication is given by
horizontal composition
\[\xymatrix{
\ar[r]^-{d_1}&.\ar[r]^-{d_0}&=\ar[r]^-{d_0d_1} &&
}\]

\[\xymatrix{&\ar@2[dd]^-{h'_{01}}&&\ar@2[dd]^-{h_{01}}
&&\\
\ar@/^1.1cm/[rr]^-{d'_1} \ar@/_1.1cm/[rr]_-{d'_0}
&&.\ar@/^1.1cm/[rr]^-{d_1}\ar@/_1.1cm/[rr]_-{d_0}&&&=\\
&&&&&}
\hskip1cm\xymatrix{&&\ar@2[dd]^-{h_{01}\,^{d_1}h'_{01}}\\
\ar@/^1.1cm/[rrrr]^-{d_1d'_1}\ar@/_1.1cm/[rrrr]_-{d_0d'_0}&&&&\\&&}
\]
etc.

Simplicial homotopy groups of $N {\mathcal C}_{(H\to D)}$ are trivial
except $\pi_0(N {\mathcal
C}_{(H\to D)}) = {\rm coker}\,\alpha$ and $\pi_1(N {\mathcal
C}_{(H\to D)}) ={\rm ker}\,\alpha$.

\subsection{Proposition} {\em $|N{\mathcal C}_{(H\to D)}|$ is the classification space
of principal $H$-bundles equipped with a chosen trivialization when the structure group is changed
to $D$ using the homomorphism $\alpha$.}
\vskip0.1cm

{\it Proof}.
In other words $|N{\mathcal C}_{(H\to D)}|$ is the homotopy fibre
of $BH \to BD$. This is the pullback under $B\alpha : BH\to BD$ of the based path
bundle $P_0 BD \to BH$ and as a principal $\Omega BD \sim D$-bundle it can
be identified with the the homotopy quotient $D//H = EH\times_\alpha D$ of
$D$ by $H$. Let us recall the $EH$ is the geometric realization of the following
simplicial space (we omit here and in all following
pictures arrows for codegeneracy maps in all following
pictures).
\[\xymatrix{&&&&&\ar[ddr]_-{h_{01}}\ar[rr]^-{h_{1}}& & \\
\ar[rr]^-{h_0} && &&
\ar@<2pt>[l]\ar@<-2pt>[l]&& &&\ldots \,\,.\\
&&&&&&\ar[uur]_{h_{0}}&}\] From here we get $EH\times_\alpha D$ as
the geometric realization of the simplicial space
\[\xymatrix{&&&&&\ar[ddr]_-{h_{01}}\ar[rr]^-{d_{1}}& & \\
\ar[rr]^-{d_0} && &&
\ar@<2pt>[l]\ar@<-2pt>[l]&& &&\ldots\\
&&&&&&\ar[uur]_{d_{0}}&}\] and we see that this really identical to the
simplicial group $N{\mathcal C}_{(H\to D)}$.

\subsection{Remark}Bundles described in proposition 4.3 are automatically left and right
$H$-principal bundles with the two principal $H$-actions commuting.
Moreover the multiplication in 4.2 gives naturally a
multiplication of such bundles. This follows from proposition 4 in
{\cite {ACJ}}. 
We will refer to such bundles as {\it crossed module
bundles}. 

If $P$ and $P'$ are two crossed module bundles and $f$ and $f'$
the corresponding classifying maps, then the point-wise product map $f.f'$ is a
classifying map for a bundle equivalent to the product bundle $P.P'$.

If $P\to X$ is a crossed module bundle then the corresponding
trivial $D$-bundle $P\times_\alpha D$ is an example of what is called a 
$(D-H)$-bundle
according to definition 5 in  {\cite {ACJ}}.

\subsection{String group}
Together with a crossed module $(H\to D)$ we can consider also
crossed modules $(H \to {\rm Im} \,\alpha)$ and $(1 \to {\rm
coker}\,\alpha)$. This gives an exact sequence of (topological)
groups
$$
1\longrightarrow|N{\mathcal C}_{(H \to {\rm Im} \,\alpha)}|
\longrightarrow |N{\mathcal C}_{(H \to D)}|
\longrightarrow|N{\mathcal C}_{(1 \to {\rm coker} \,\alpha)}|= {\rm
coker} \,\alpha \longrightarrow 1$$ String group {\it String} is a
nice example of the above construction. Let $G$ be a simply
connected compact simple Lie group. The crossed module in question
is $H=\widehat{\Omega G}$ the centrally extended group of based
loops and $D=P_0G$ is the group of based paths \cite{BCSS}, \cite{H}
or some modification of these \cite{ST}.

Of course we can as well consider crossed modules $({\ker \alpha}
\to e)$ and $({\rm Im}\,\alpha)\to D)$ in which case we obtain the
exact sequence
$$
1\longrightarrow|N{\mathcal C}_{({\ker \alpha} \to e)}|=B\ker \alpha
\longrightarrow |N{\mathcal C}_{(H \to D)}|
\longrightarrow|N{\mathcal C}_{({\rm Im}\,\alpha\to D)}|
\longrightarrow 1.$$ We can view homotopy quotient $D//H =
EH\times_\alpha D$ as a bundle with the base space ${\rm coker}
\,\alpha$. $EH$ is the universal bundle for any subgroup of $H$ and hence for
the normal subgroup $\ker \alpha$ too. The action of $H$ on $EH$ descents to an action of
$H/\ker \alpha \sim \alpha(H)$ on $B\ker \alpha$ and we see that we have the
bundle $B\ker \alpha \times_{\alpha(H)} D$. The two
exact sequences in  remark 4.4 are identical to
$$
1\longrightarrow B\ker \alpha \longrightarrow |N{\mathcal C}_{(H \to
D)}| \longrightarrow {\rm coker} \,\alpha \longrightarrow 1.$$

There is another nice description of the group structure on 
$EH\times_\alpha D$\footnote{I thank D. Stevenson
for noticing this to me}.
$EH$ itself can be thought of as $|N{\mathcal C}_{(H \to H)}|$. Hence it is
a topological group. The action of $D$ on $H$ naturally extends to $EH$
and we can form the semidirect product $EH \rtimes D$. This group structure 
factors to $EH \rtimes_\alpha D$. 

Now if we equip $B\ker \alpha$ with the factor group structure then the
$D$-action factors to $B\ker \alpha$ as it preserves $\ker \alpha$.
It is easy to check that  $EH \rtimes_\alpha D$ and 
$B\ker \alpha \rtimes_{\alpha(H)} D$ are identical as topological groups.

This description of {\em String} is very close to the one of Stolz and
Teichner \cite{ST}. Very briefly, in their construction (of {\em String})
$H$ is $\tilde L_IG$, the central extension of $L_I G$ (group
of all piece-wise smooth loops $\gamma: S^1 \to G$
with the support in the upper semicircle $I \in S^1$). Here $G$ is a compact,
simply connected Lie group. Their $D$ is the
group of based paths
$P_e^I G=\{\gamma : I \to G| \gamma (1)=e\}$. With these choices they can take 
$PU(A_\rho)$ as
a model for $B\ker \alpha$, where $A_\rho$ is certain von Neumann
algebra (type ${\rm III}_1$ factor) associated with the vacuum representation of
the loop group $LG$ at some fixed level $l \in H^4(BG)$. See \cite{ST} for
details. Related discussion in terms of Morita equivalence of 2-groups will
appear in \cite{danny}.

\subsection{Crossed module as a 2-category}
Similarly we denote $\tilde {\mathcal C}_{(H \to D)}$ the
(topological) 2-category with just one object, 1-arrows group
elements $d \in D$ and 2-arrows group elements $(h,d)$ of $H\rtimes
D $. Again we can form the corresponding nerve $N\tilde {\mathcal
C}_{(H\to D)}$ \cite{D}. This simplicial manifold can be pictorially
represented as:
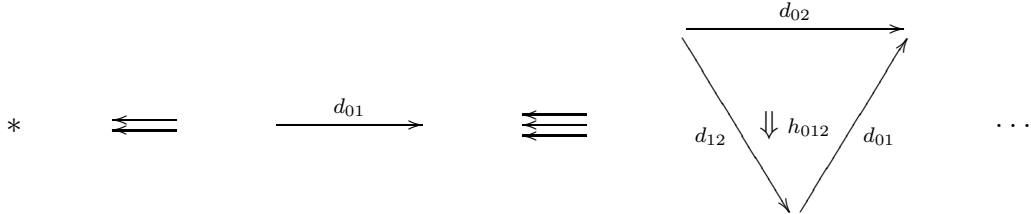
\[\xymatrix{&&&&&&&&\ar[ddr]_-{d_{12}}\ar[rr]^-{d_{02}}& & \\
\ast &&\ar@<2pt>[l]\ar@<-2pt>[l]& \ar[rr]^-{d_{01}} && &&
\ar@<4pt>[l]\ar@<-4pt>[l]\ar[l]&&\Downarrow{\scriptstyle h_{012}} &&\ldots\\
&&&&&&&&&\ar[uur]_{d_{01}}&}\]
Simplicial homotopy groups of $N\tilde {\mathcal C}_{(H\to D)}$ are trivial
except $\pi_1(N\tilde {\mathcal
C}_{(H\to D)}) = {\rm coker}\,\alpha$ and $\pi_2(N\tilde {\mathcal
C}_{(H\to D)}) ={\rm ker}\,\alpha$.

\section{Bundle gerbes}
Let us recall
the definition of a nonabelian bundle gerbe in the form given in \cite{ACJ}.
Consider a submersion $\wp~:~Y\rightarrow X$ (i.e. a map onto with
differential onto). It follows we can always find an open covering $\{{O}_{\alpha}\}$ of $M$
with local sections $\sigma_{\alpha}:{O}_\alpha\rightarrow Y$,
i.e. $\wp\circ \sigma_\alpha=id$. The manifold $Y$ will always be 
equipped with the submersion $\wp~:~Y\rightarrow X$. 
We also consider $Y^{[n]}=Y\times_X Y\times_X Y\ldots\times_X Y$ the n-fold
fiber product of $Y$, i.e. 
$Y^{[n]}\equiv\{(y_1,\ldots y_n)\in Y^n\;|\;\wp(y_1)=\wp(y_2)=
\ldots\wp(y_n)\}$. 

Given a  $(H\to D)$-crossed module bundle $\mathcal E$ over $Y^{[2]}$
we denote by ${\mathcal E}_{12}=p_{12}^*(\mathcal E)$ 
the crossed module bundle on $Y^{[3]}$ obtained as pull-back 
of $p_{12}:Y^{[3]}\rightarrow Y^{[2]}$ ($p_{12}$ is the identity on
its first two arguments);
similarly for ${\mathcal E}_{13}$ and ${\mathcal E}_{23}$.

Consider the quadruple
$({\mathcal E},Y,X,\hh)$ where ${\mathcal E}$ is a crossed module bundle, $Y\to M$ a submersion
and $\hh$ an isomorphism of
crossed module bundles
$\hh :  {\mathcal E}_{12}{\mathcal E}_{23}\to{\mathcal E}_{13}$ (let us recall that two crossed module
bundles can be multiplied to obtain again a crossed module bundle). 
We now consider 
$Y^{[4]}$ and the bundles ${\mathcal E}_{12},\,{\mathcal E}_{23},\,{\mathcal E}_{13},\,{\mathcal E}_{24},
\,{\mathcal E}_{34},\,{\mathcal E}_{14}$
on $Y^{[4]}$ 
relative to the projections $p_{12}~:~Y^{[4]}\rightarrow Y^{[2]}$ etc.
and also the crossed module isomorphisms
$\hh_{123}, \,\hh_{124},\, \hh_{123}, \,\hh_{234}$ induced by projections 
$p_{123}:Y^{[4]}\rightarrow Y^{[3]}$ etc.

\subsection{Definition} 
The quadruple
$({\mathcal E},Y,X,\hh)$ is called a crossed module bundle gerbe if $\hh$ satisfies the
cocycle condition (associativity) on $Y^{[4]}$
\begin{center}
$$
\begin{CD}
 \EE_{12}\EE_{23}\EE_{34}@> 
 {\hh_{234}}>> \EE_{12}\EE_{24}\\
 @V {\hh_{123}}
 VV @V  VV \hskip-0.3cm{\scriptstyle{\hh_{124}}}\\ 
\EE_{13}\EE_{34}@> 
 {\hh_{134}}>>\EE_{14}\,.
\end{CD}
$$
\end{center}

\subsection{Definition}
Two crossed module bundle gerbes $(\EE,Y,X,\hh)$ and $(\EE',Y',X,\hh')$ are
stably isomorphic if there exist a crossed module bundle $\NN \to Z=Y\times_X
Y'$ such that over $Z^{[2]}$ the crossed module bundles $q^*\EE \NN_2$ and 
$\NN_1 q'^*\EE'$ are isomorphic. The corresponding isomorphism
$\ll: q^*\EE \NN_2 \to \NN_1 q'^*\EE'$ should satisfy on $Y^{[3]}$ the condition
$$\ll_{13} \circ \hh = \hh' \circ \ll_{12} \circ \ll_{23}.$$
Here $q$ and $q'$ are projections onto first and second factor of $Z=Y\times_X
Y'$. $\NN_1$ and $\NN_2$ are the pullbacks of  $\NN \to Z$ to $Z^{[2]}$ under respective
projections form $Z^{[2]}$ to $Z$ etc.

\subsection{Remark}
Locally bundle gerbes can be described in terms of cocycles as follows.
We can consider the trivializing cover $\{O_\alpha \}$ of the submersion
$Y\to X$ be a good one.
Then a crossed module bundle gerbe
can be described by a cocycle $\{d_{\alpha \beta}, h_{\alpha \beta
\gamma}\}$ where the maps $d_{\alpha \beta}: O_\alpha \cap O_\beta
\to D$ and $h_{\alpha \beta \gamma}: O_\alpha \cap O_\beta \cap
O_\gamma \to H$ fulfill the following cocycle condition
$$d_{\alpha \beta}d_{\beta \gamma}= \alpha(h_{\alpha \beta \gamma})
d_{\alpha \gamma} \hskip0.4cm {\rm on}\hskip0.4cm O_\alpha \cap O_\beta \cap O_\gamma
$$
and
$$h_{\alpha \beta \gamma}h_{\alpha \gamma \delta} = \,^{d_{\alpha
\beta}}\hskip-0.05cm
h_{\beta \gamma \delta} h_{\alpha \beta \delta} \hskip0.4cm {\rm on}\hskip0.4cm
O_\alpha \cap O_\beta \cap O_\gamma \cap O_\delta.$$

Two crossed module bundle gerbes are stably equivalent if
their respective cocycles
 $\{d_{\alpha \beta}, h_{\alpha \beta \gamma}\}$ and
 $\{d_{\alpha \beta}', h_{\alpha \beta \gamma}'\}$ are related
 as
 $$d_{\alpha \beta}'=d_{\alpha}\alpha(h_{\alpha \beta})
 d_{\alpha \beta}d_{\beta}^{-1}$$
 and
 $$ h_{\alpha \beta \gamma}'= \,^{d_\alpha}h_{\alpha \beta}
 \,^{d_\alpha d_{\alpha \beta}}h_{\beta \gamma}
 \,^{d_\alpha}h_{\alpha \beta \gamma}\,^{d_\alpha}h_{\alpha \beta}^{-1}$$
 with $d_\alpha : O_\alpha \to D$ and $h_{\alpha\beta} : O_\alpha \cap O_\beta
 \to H$.

Pullback of a bundle gerbe is obtained pulling back
the corresponding cocycle.

\subsection{Universal $N{\mathcal C}_{(H\to D)}$ bundle}
In 4.2 we have described the simplicial group $N{\mathcal
C}_{(H\to D)}$. Now we can construct the corresponding universal
bundle. As a result we get simplicial manifolds $\overline
WN{\mathcal C_{H\to D}}$ and $WN{\mathcal C_{H\to D}}$ which are
pictorially represented as
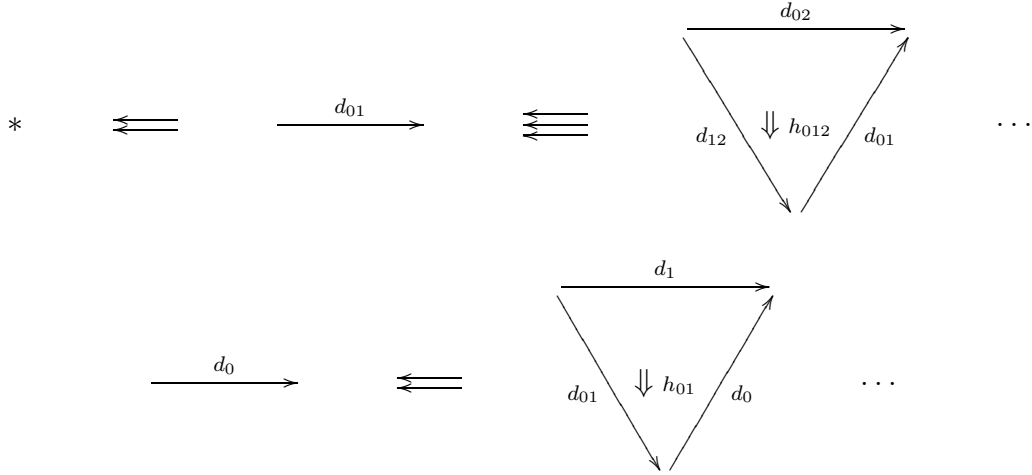
\[\xymatrix{&&&&&&&&\ar[ddr]_-{d_{12}}\ar[rr]^-{d_{02}}& & \\
\ast &&\ar@<2pt>[l]\ar@<-2pt>[l]& \ar[rr]^-{d_{01}} && &&
\ar@<4pt>[l]\ar@<-4pt>[l]\ar[l]&&\Downarrow{\scriptstyle h_{012}} &&\ldots\\
&&&&&&&&&\ar[uur]_{d_{01}}&}\]
\[\xymatrix{&&&&&\ar[ddr]_-{d_{01}}\ar[rr]^-{d_{1}}& & \\
\ar[rr]^-{d_0} && &&
\ar@<2pt>[l]\ar@<-2pt>[l]&&\Downarrow{\scriptstyle h_{01}} &&\ldots\\
&&&&&&\ar[uur]_{d_{0}}&}\]

Comparing to 4.6 gives
\subsection{Proposition}
{\em $\overline WN{\mathcal C}_{(H\to D)}= N\tilde {\mathcal C}_{(H\to D)}$.}
\vskip0.1cm
Now we can touch upon the question of the classification of crossed module
bundle gerbes.
\subsection{Theorem} {\em
Equivalence classes of  principal $|N{\mathcal C}_{(H\to
D)}|$-bundles are one to one with stable equivalence classes of
$(H\to D)$ crossed module bundle gerbes. The geometric realization $|WN
{\mathcal C}_{(H\to D)}|=E|N{\mathcal C}_{(H\to D)}| \to |\overline WN
{\mathcal C}_{(H\to D)}|=B|N{\mathcal C}_{(H\to D)}|$ gives the
universal $|N{\mathcal C}_{(H\to D)}|$-bundle as well as the
universal crossed module bundle gerbe. }
\vskip0.1cm

{\it{Proof.}}
The proof is just a slight generalization of section 5 in \cite{ACJ}, where the
lifting bundle gerbes (crossed module bundle gerbes with $\ker \alpha =0$) are 
discussed in detail.   
Let $f:X \to B|N{\mathcal C}_{(H\to D)}|$ be the classification map for an
$|N{\mathcal C}_{(H\to D)}|$-principal bundle $P$. 
Associated with $P$ there is a map $P^{[2]}\to |N{\mathcal C}_{(H\to D)}|$
which sends $(p,p')\in P$ in the same fibre into unique group element
$g\in |N{\mathcal C}_{(H\to D)}|$ which relates $p$ and $p'$.
As $|N{\mathcal C}_{(H\to D)}|$ is the classification space for crossed module
bundles, we obtain that way a crossed module bundle $\EE \to P^{[2]}$.
As it follow from   remark 4.4, $\EE_{12}\EE_{23}$ is isomorphic to $\EE_{13}$, and
it is easy to check that this isomorphism fulfils the cocycle condition of
definition 5.1.
So we obtain a bundle gerbe with $Y=P$.
If we start with an equivalent bundle $P'$ we obtain an stably equivalent gerbe.

Conversely if we start with a crossed module bundle gerbe, the classification
map of crossed module bundle $\EE \to Y^{[2]}$ is a map from 
$f:Y^{[2]} \to 
|N{\mathcal C}_{(H\to D)}|$ fulfilling the on $Y^{[3]}$ the cocycle condition
$f(y_1, y_2)f(y_2,y_3)=f(y_1,y_3)$. Using local sections $O_\alpha\to Y$ we get
a cocycle
$f_{\alpha \beta}: O_\alpha \cap O_\beta \to |N{\mathcal C}_{(H\to D)}|$; 
$f_{\alpha \beta}f_{\beta \gamma}=f_{\alpha \gamma}$. Thus
we have transition functions for an $|N{\mathcal C}_{(H\to D)}|$-bundle. Starting form a
stably equivalent gerbe we get an equivalent bundle.

\subsection{Remark} Let us recall that by definition under a nonabelian
$H$-bundle gerbe we understand an $(H\to {\rm
Aut}(H))$-crossed module bundle gerbe \cite{ACJ}. So the universal $H$-bundle
gerbe is the same as the universal $|N{\mathcal C}_{(H\to {\rm
Aut}(H))}|$-bundle.

\subsection{String structures}
Now can apply the classifying space functor $B$ to the exact
sequence of 4.5, or which is the same the $\overline W$ to the
corresponding simplicial groups.
Hence have the following exact sequence ($\ker \alpha$ is abelian)
$$
1\to B{\rm
ker}\,\alpha\to
|N{\mathcal C}_{(H \to D)}|
 \to {\rm coker} \,\alpha  \to
 B^2{\rm
ker}\,\alpha \to
B|N{\mathcal C}_{(H \to D)}|
 \to B{\rm coker} \,\alpha \to B^3{\rm
ker}\,\alpha.
$$
It follows that a lift of a principal ${\rm coker}\,\alpha$-bundle
to a principal $|N{\mathcal C}_{(H\to D)}|$-bundle is the same as a
lift of an $(\alpha (H)\to D)$-bundle gerbe to an $(H \to D)$-bundle
gerbe.

In the case of {\it String} we do have
$$
1\to K(\BZ,2) \to String \to Spin \to K(\BZ,3) \to BString \to BSpin \to
K(\BZ,4).
$$

String structure is a lift of the structure group of a principal
$Spin$-bundle to the string group $String$ \cite{ST}. So the string structure
is also lift of an $(\Omega Spin \to P_0
Spin)$-bundle gerbe  to an  $(\widehat{\Omega Spin} \to P_0
Spin)$-bundle gerbe.

\subsection{Remark} A crossed module bundle gerbe is canonically equipped with
a module (see section 6 of \cite{ACJ} for the definition of a bundle gerbe
module). The trivial $D$-principal bundle $D\times Y \to Y$ fulfils all the
axioms of a module. This is shown in \cite{ACJ} in the case $D={\rm Aut}(H)$
and applies word by word to the more general situation as well.

\subsection{Remark}
Let us consider  the (topological) 1-category (actually 1-groupoid) 
${\mathcal C}_{\{O_\alpha\}}$, related to an open covering ${\{O_\alpha\}}$,
described as follows. Objects are
pairs $(x, O_\alpha)$
with $x\in O_\alpha$ and there is unique morphism $(x, O_\alpha) \to (y,
O_\beta)$ iff $x=y \in O_\alpha \cap O_\beta$.   Let 
$N{\mathcal C}_{\{O_\alpha\}}$ denote the
nerve of this category. Consider a map 
of simplicial sets $ N{\mathcal C}_{\{O_\alpha\}}\to \overline WN {\mathcal
C}_{(H\to D)}$.  Then the maps between 1- 2- and 3-simplexes
 give us the cocycle  for gerbe transition
functions in the definition 5.1. We also see that the simplicial $\tau_1$  is identified with
$d_{\alpha
\beta}$, $\tau_2$ identifies with $d_{\alpha \gamma}d_{\beta
\gamma}^{-1} \stackrel{h_{\alpha \beta \gamma}}{\longrightarrow}
d_{\alpha \beta}$ etc.  So from 1.3 and 5.3 we can conclude that locally the stable 
equivalence classes of
crossed module gerbes are described by homotopy classes of simplicial maps
$N{\mathcal C}_{\{O_\alpha\}}\to \overline{W}N{\mathcal C}_{(H\to D)}= 
N\tilde{\mathcal C}_{(H\to D)}$.

\section{Connection and $B$-field on a bundle gerbe} In the previous section we have established a
correspondence between $|N{\mathcal C}_{(H\to D)}|$-principal
bundles and $(H\to D)$-crossed module bundle gerbes. Now we
would like to extend this relationship to connections, and also
discuss the $B$-field from this point of view. However let us recall
that $|N{\mathcal C}_{(H\to D)}|$ is only a topological group so in
general there is no differential geometric connection on a
principal $|N{\mathcal C}_{(H\to D)}|$-bundle over a manifold $X$.
But we have the simplicial connection as described in section 2 on
any simplicial $N{\mathcal C}_{(H\to D)}$-bundle $P\to X$.

The notion of a bundle gerbe connection (and that of a bundle gerbe $B$-field
as well) are
quite subtle and we are not going to repeat them here in their global
formulations (see \cite{ACJ}, \cite{BM} for that). Instead we will give their 
local description using cocycles.
This description is perfectly well suited for our purposes as we will relate
the bundle gerbe connection and $B$-field to the simplicial connection and
simplicial $\bar B$-field as they were introduced in sections 2 and 3
in the case of a simplicial $N{\mathcal C}_{(H\to D)}$-bundle
over $N{\mathcal C}_{\{O_\alpha\}}$ described by a classifying map 
$N{\mathcal C}_{\{O_\alpha\}}\to \overline{W}N{\mathcal C}_{(H\to D)}= 
N\tilde{\mathcal C}_{(H\to D)}$ (see remark 5.10).

Let us now recall the cocycle description of what a connection on
an crossed module bundle gerbe is. Again let $\{O_{\alpha}\}$
be an open covering of a manifold $X$.
\subsection{Bundle gerbe connection}
A collection $\{A_\alpha, a_{\alpha \beta}\}$, with $A_\alpha \in
\Omega^1(O_\alpha)\otimes {\rm Lie}(D)$ and $a_{\alpha \beta}\in
\Omega^1(O_\alpha \cap O_\beta)\otimes {\rm Lie}(H) $ is called a
connection on crossed module bundle gerbe (characterized by a nonabelian cocycle $\{d_{\alpha\beta}, h_{\alpha
\beta \gamma}\}$) if it fulfils the following conditions
$$A_\alpha = d_{\alpha \beta}A_\beta  d_{\alpha \beta}^{-1} +
d_{\alpha \beta} d d_{\alpha \beta}^{-1} + \alpha(a_{\alpha \beta})
\hskip0.4cm {\rm on} \hskip0.4cm O_\alpha \cap O_\beta$$
and
$$a_{\alpha \beta} + ^{d_{\alpha \beta}}\hskip-0.1cm a_{\beta \gamma}=
h_{\alpha \beta \gamma}a_{\alpha \gamma}h_{\alpha \beta \gamma}^{-1}
+ h_{\alpha \beta \gamma}dh_{\alpha \beta \gamma}^{-1} + 
T_{A_\alpha}(h_{\alpha \beta
\gamma}^{-1})\hskip0.4cm {\rm on}\hskip0.4cm O_\alpha \cap O_\beta
\cap O_\gamma.$$
Here for $A$ a ${\rm Lie}(D)$-valued one form and $h\in H$ the  
${\rm Lie}(H)$-valued one form $T_A(h)$ is defined as follows. For 
$X \in {\rm Lie}(D)$ we put $T_X (h)= [h\,^{{\rm exp}(tX)}(h^{-1})]$, where
the bracket $[\,\,\,\,]$ means the tangent vector to the curve at the group
identity $1_H$. For ${\rm Lie}(D)$-valued one form $A=A^\rho X^\rho$, with
$\{X^\rho\}$ a basis of  ${\rm Lie}(D)$, we put $T_A \equiv A^\rho T_{X^\rho}$.

The curvature $F$ is given by a collection of local two-forms
$F_\alpha \in \Omega^2(O_\alpha)\otimes {\rm Lie}(D)$ defined as
$F_\alpha = dA_\alpha + A_\alpha \wedge A_\alpha$; the corresponding
cocycle conditions follow from the definition. We will not repeat
the explicit formulas here, interested reader can find them in e.g.
\cite{ACJ}, \cite{BM}. Now we can compare the above definition with the definition
of a simplicial connection on a $N{\mathcal C}_{(H\to D)}$-principal
bundle $P\to N{\mathcal C}_{\{O_\alpha\}}$. Realizing that $\tau_1$ corresponds
$d_{\alpha \beta}$, $\tau_2$ corresponds to $d_{\alpha
\gamma}d_{\beta \gamma}^{-1} \stackrel{h_{\alpha \beta
\gamma}}{\longrightarrow} d_{\alpha \beta}$, $A_0$ corresponds to
$A_\alpha$, $a_{01}$ of $A_1=(\partial_0 A_1 \stackrel{a_{01}}{\to} 
\partial_1 A_1)$ corresponds to $-a_{\alpha \beta}$ etc., we easily obtain
\subsection{Proposition}
{\em A connection on a crossed module bundle gerbe
defines a simplicial connection on the corresponding $N{\mathcal
C}_{(H\to D)}$-principal bundle over $N{\mathcal C}_{\{O_\alpha\}}$ 
and vice versa.}
\vskip.1cm
Similar discussion applies to $B$-field as well.
\subsection{Bundle gerbe $B$-field} $B$-field on a crossed module bundle gerbe
equipped with a connection is a collection $\{B_\alpha,
\delta_{\alpha \beta}\}$ of local two-forms $B_\alpha \in
\Omega^2(O_\alpha) \otimes {\rm Lie}(H)$ and $\delta_{\alpha \beta}\in
\Omega^2(O_{\alpha \beta})\otimes {\rm Lie}(H)$ such that
$$ B_\alpha = \,^{d_{\alpha \beta}}B_\beta + \delta_{\alpha \beta}\hskip0.4cm {\rm on} \hskip0.4cm
O_{\alpha}\cap O_\beta$$ and
$$\delta_{\alpha \beta} +\,^{d_{\alpha \beta}}\delta_{\beta \gamma}=
h_{\alpha \beta \gamma} \delta_{\alpha \gamma}h_{\alpha \beta
\gamma}^{-1}  +B_\alpha -h_{\alpha \beta \gamma}B_\alpha h_{\alpha
\beta \gamma}^{-1}\hskip0.4cm {\rm on} \hskip0.4cm O_{\alpha}\cap
O_\beta\cap O_\gamma.$$

Given a simplicial $\bar B$ in the present case then the bundle
gerbe $B$-field is identified as the morphism $B$ in the $\bar B_1=
(\partial_0 \bar B_1 
\stackrel{- B}{\longrightarrow}
0)$ part of the simplicial $\bar B$ and the simplicial $(\partial_2 \bar B_2 -
\partial_1 \bar B_2)$
is identified with the bundle gerbe $\delta$, we obtain the
following proposition.
\subsection{Proposition}
{\em A simplicial $\bar B$-field on a $N{\mathcal C}_{(H\to
D)}$ principal bundle over $N{\mathcal C}_{\{O_\alpha\}}$ gives a $B$-field on the 
corresponding bundle gerbe and vice versa. }
\vskip0.1cm
The bundle gerbe $\nu$-field is  defined as $\nu = F + \alpha(B)$.
This definition guaranties that it is the same as the simplicial one
in the present case.

\subsection{Remark}
It is generally true only in the case of abelian $H$ that
connection $A$ and the $B$-field can be chosen such that $\nu_\alpha =
d_{\alpha \beta}\nu_\beta d_{\alpha \beta}^{-1}$. We are are not sure what kind
of condition should replace this in the case of nonabelian $H$.

\subsection*{Acknowledgement}
I would like to thank P. Aschieri, I. Bakovi\'c, J. Baez, D. Husem\"oller, D.
Roberts, U. Schreiber, and especially A. Henriques and D.
Stevenson for interesting and helpful discussions.

\end{document}